\documentclass[12pt]{article}
\usepackage{bbm}
\usepackage{latexsym}
\usepackage{mathrsfs}
\usepackage{graphicx}
\usepackage{subfigure}
\usepackage{amssymb}
\usepackage{amsmath}
\usepackage{amsthm}
\usepackage{color}
\usepackage{cite}
\usepackage{float}
\usepackage{indentfirst}
\usepackage{anysize}\marginsize{45mm}{45mm}{40mm}{50mm}
\usepackage{caption}
\usepackage{float}
\usepackage{multirow}

\usepackage{enumerate}

\usepackage{latexsym, bm}
\newtheoremstyle{lemma}{\topsep}{\topsep}%
     {}%         Body font
     {}%         Indent amount (empty = no indent, \parindent = para indent)
     {\bfseries}% Thm head font
     {}%        Punctuation after thm head
     {0.1em}%     Space after thm head (\newline = linebreak)
     {\thmname{#1}\thmnumber{ #2}\thmnote{ #3}}%         Thm head spec
\theoremstyle{lemma}  %%定理引用参考文献格式

\newtheorem{theorem}{Theorem}     %全局定义
\newtheorem{lemma}[theorem]{Lemma}

\newtheorem{definition}{Definition}

\numberwithin{equation}{section}

\title{ Structure and substructure connectivity of balanced hypercubes\thanks{This research was partially supported by the National Natural Science Foundation of China (No. 11761056), the Chunhui Project of Ministry of Education (No. Z2017047) and the fundamental research funds for the central universities (No. 2672018ZYGX2018J069)}}

\author{ Huazhong L\"{u}$^{1}$\thanks{Corresponding author.} and Tingzeng Wu$^{2}$\\
{\small $^{1}$School of Mathematical Sciences, University of Electronic Science and Technology of China,} \\
{\small Chengdu, Sichuan 610054, P.R. China}\\
{\small E-mail: lvhz08@lzu.edu.cn}\\
{\small $^{2}$School of Mathematics and Statistics, Qinghai Nationalities University, }\\
{\small Xining, Qinghai 810007, P.R. China} \\
{\small E-mail: mathtzwu@163.com}\\}
\date{}
\begin{document}

\maketitle
\begin{abstract}

The connectivity of a network directly signifies its reliability and fault-tolerance. Structure and substructure connectivity are two novel generalizations of the connectivity.
Let $H$ be a subgraph of a connected graph $G$. The structure connectivity (resp. substructure connectivity) of $G$, denoted by $\kappa(G;H)$ (resp. $\kappa^s(G;H)$), is defined to be the minimum cardinality of a set $F$ of connected subgraphs in $G$, if exists, whose removal disconnects $G$ and each element of $F$ is isomorphic to $H$ (resp. a subgraph of $H$). In this paper, we shall establish both $\kappa(BH_n;H)$ and $\kappa^s(BH_n;H)$ of the balanced hypercube $BH_n$ for $H\in\{K_1,K_{1,1},K_{1,2},K_{1,3},C_4\}$.

\noindent \textbf{Key words:} Interconnection networks; Structure connectivity; Substructure connectivity; Balanced hypercube

\noindent \textbf{Mathematics Subject Classification:} 05C40, 68R10
\end{abstract}

\section{Introduction}

The interconnection network is crucial in parallel processing and distributed system since the performance of the system is significantly determined by its topology. As the size of a network increases continuously, the reliability and fault-tolerance become central issues. The classical connectivity is an important measure to evaluate fault-tolerance of a network with few processors. An obvious deficiency of the connectivity is the assumption that all the parts of the network can be potentially fail at the same time. However, in large networks, it is unlikely that all the vertices incident to a vertex fail simultaneously, indicating high resilience of large networks. To address the shortcomings of the connectivity stated above, Harary \cite{Harary} introduced the conditional connectivity of a connected graph by adding some constraints on the components of the resulting graph after vertex deletion. After that, several kinds of conditional connectivity were proposed and investigated \cite{Chen,Esfahanian,Fabrega,Hsieh,Ning2,Wan,Ye}, such as $g$-connectivity and $h$-connectivity.

The $g$-connectivity of $G$, denoted by $\kappa_0^g(G)$, if exists, is defined as the minimum cardinality of a vertex set in $G$, if exists, whose deletion disconnects $G$ and leaves each remaining component with at least $g+1$ vertices. The $h$-connectivity of $G$, denoted by $\kappa^h(G)$), if exists, is defined as the cardinality of a minimum cardinality of a vertex set in $G$, if exists, whose deletion disconnects $G$ and each vertex in the resulting graph has at least $h$ neighbors. From the definitions above, it is obvious that $\kappa_0^g(G)\leq\kappa_0^{g+1}(G)$ and $\kappa^h(G)\leq \kappa^{h+1}(G)$ if $G$ has $\kappa_0^{g+1}(G)$ and $\kappa^{h+1}(G)$, respectively. So both of $g$-connectivity and $h$-connectivity are generalizations of the connectivity, which supply more accurate measures to evaluate reliability and fault-tolerance of large networks. Moreover, the higher $g$-connectivity or $h$-connectivity the network has, the more reliable the network is \cite{Esfahanian,Lu2}. It is known that there exists no polynomial time algorithm to compute the $g$-connectivity and $h$-connectivity of a general graph \cite{Chang,Esfahanian}. The $h$-connectivity \cite{Chen,Li,Ning2,Oh,Wan,Ye} and $g$-connectivity \cite{Chen.L,Hsieh,Li2,L.Xu,W.Yang,Zhang} of some famous networks are investigated in the literature.

As stated above, most studies on reliability and fault-tolerance of networks are under the assumption that the status of a vertex $u$, whether it is good or faulty, is an event independent of the status of vertices around $u$. In other words, vertices that are linked in a network do not affect each other. Nevertheless, in reality, the neighbors of a faulty vertex might be more vulnerable or have a higher possibility of becoming faulty later. Also note that networks and subnetworks are made into chips. This means that when any vertex is faulty, the whole chip is regarded as faulty. Motivated by these, Lin et al. \cite{Lin} proposed structure and substructure connectivity to evaluate the fault-tolerance of networks not only from the perspective of individual vertex, but also some special structure of the network.

A set $F$ of connected subgraphs of $G$ is a subgraph-cut of $G$ if $G-V(F)$ is disconnected or trivial. Let $H$ be a connected subgraph of $G$, then $F$ is an $H$-structure-cut if $F$ is a subgraph-cut, and each element in $F$ is isomorphic to $H$. The $H$-structure-connectivity of $G$, denoted by $\kappa(G; H)$, is the minimum cardinality of all $H$-structure-cuts of $G$. Furthermore, $F$ is an $H$-substructure-cut if $F$ is a subgraph-cut, such that each element in $F$ is isomorphic to a connected subgraph of $H$. The $H$-substructure-connectivity of $G$, denoted by $\kappa^s(G; H)$, is the minimum cardinality of all $H$-substructure-cuts of $G$.

The balanced hypercube was proposed by Wu and Huang \cite{Wu} as a novel interconnection network. As an alternative of the well-known hypercube, the balanced hypercube keeps lots of desirable properties of the hypercube, such as bipartite, high symmetry, scalability, etc. It is known that odd-dimension balanced hypercube has a smaller diameter than that of the hypercube of the same order. In particular, the balanced hypercube is superior to the hypercube in a sense that it supports an efficient reconfiguration without changing the adjacent relationship among tasks \cite{Wu}. Some other excellent properties of the balanced hypercube were discussed by many researchers, such as fault-tolerant resource placement problem \cite{Huang2} $g$-connectivity \cite{Yang2,Lu2,D-Yang} and $h$-connectivity \cite{Lu4}, Hamiltonian path (cycle) embedding \cite{Cheng,Hao,Xu,Yang,P.Li}, matching preclusion \cite{Lu} and matching extendability \cite{Lu3}, conditional diagnosability \cite{Yang3} and symmetric properties \cite{Zhou,Zhou2}.

\vskip 0.05 in

Lin et al. \cite{Lin} considered $\kappa(Q_n;H)$ and $\kappa^s(Q_n;H)$ of the hypercube $Q_n$ for $H\in\{K_1, K_{1,1}, K_{1,2}, K_{1,3},C_4\}$. Later, Sabir and Meng \cite{Sabir} generalized the results in $Q_n$ and studied this problem in the folded hypercube. Mane \cite{Mane} determined $\kappa(Q_n; Q_m)$ with $m\leq n-2$ and obtained the upper bound of $\kappa(Q_n; C_{2k})$ with $2\leq k\leq 2^{n-1}$. Furthermore, Lv et al. \cite{Lv} investigated $\kappa(Q_n^k;H)$ and $\kappa^s(Q_n^k;H)$ of the $k$-ary $n$-cube hypercube $Q_n^k$ for $H\in\{K_1, K_{1,1}, K_{1,2}, K_{1,3}\}$. In this paper, we will establish $\kappa(BH_n;H)$ and $\kappa^s(BH_n;H)$ of the balanced hypercube $BH_n$ ($n\geq2$) for $H\in\{K_1, K_{1,1}, K_{1,2}, K_{1,3}, C_4\}$. Note that $K_1$ is a singleton, the $K_1$-structure connectivity degenerate to traditional connectivity.

\vskip 0.05 in

The rest of this paper is organized as follows. In Section 2, the definitions of balanced hypercubes and some useful lemmas are presented. The main results of this paper are shown in Section 3. Conclusions are given in Section 4.

\vskip 0.05 in

\section{Preliminaries}

Let $G=(V(G),E(G))$ be a graph, where $V(G)$ is vertex-set of $G$ and $E(G)$ is edge-set of $G$. The number of vertices of $G$ is denoted by $|G|$. The {\em neighborhood} of a vertex $v$ is the set of vertices adjacent to $v$, written as $N_{G}(v)$. Let $F\subseteq V(G)$, we define $N_{G}(F)=\cup_{v\in F}N_{G}(v)-F$. For $A\subset G$, we use $N_{G}(A)$ to denote $N_{G}(V(A))$ briefly. For other standard graph notations not defined here please refer to \cite{Bondy}.

In what follows, we shall give definitions of the balanced hypercube and some lemmas.
\vskip 0.0 in

\begin{definition}{\bf .}\label{def1}\cite{Wu} An $n$-dimensional balanced hypercube $BH_{n}$ consists of
$2^{2n}$ vertices $(a_{0},\ldots,a_{i-1},a_{i},a_{i+1},\ldots,a_{n-1})$, where $a_{i}\in\{0,1,2,3\}(0\leq
i\leq n-1)$. An arbitrary vertex $v=(a_{0},\ldots,a_{i-1},$
$a_{i},a_{i+1},\ldots,a_{n-1})$ in $BH_{n}$ has the following $2n$ neighbors:

\begin{enumerate}
\item $((a_{0}+1)$ mod $
4,a_{1},\ldots,a_{i-1},a_{i},a_{i+1},\ldots,a_{n-1})$,\\
      $((a_{0}-1)$ mod $ 4,a_{1},\ldots,a_{i-1},a_{i},a_{i+1},\ldots,a_{n-1})$, and
\item $((a_{0}+1)$ mod $ 4,a_{1},\ldots,a_{i-1},(a_{i}+(-1)^{a_{0}})$ mod $
4,a_{i+1},\ldots,a_{n-1})$,\\
      $((a_{0}-1)$ mod $ 4,a_{1},\ldots,a_{i-1},(a_{i}+(-1)^{a_{0}})$ mod $
      4,a_{i+1},\ldots,a_{n-1})$.
\end{enumerate}
\end{definition}

The first coordinate $a_{0}$ of the vertex
$(a_{0},\ldots,a_{i},\ldots,a_{n-1})$ in $BH_{n}$ is defined as {\em inner index}, and
other coordinates $a_{i}$ $(1\leq i\leq n-1)$ {\em outer index}.

\vskip 0.0 in

The following definition shows recursive property of the balanced hypercube.

\begin{definition}{\bf .}\label{def2}\cite{Wu}
\begin{enumerate}
\item $BH_{1}$ is a $4$-cycle and the vertices are labelled
by $0,1,2,3$ clockwise.
\item $BH_{k+1}$ is constructed from four $BH_{k}$s, which
are labelled by $BH^{0}_{k}$, $BH^{1}_{k}$, $BH^{2}_{k}$,
$BH^{3}_{k}$. For any vertex in $BH_{k}^{i}(0\leq i\leq 3)$, its new labelling in $BH_{k+1}$ is $(a_{0},a_{1},\ldots,a_{k-1},i)$, and it has two new neighbors:
\begin{enumerate}
\item[a)] $BH^{i+1}_{k}:((a_{0}+1)$mod $4,a_{1},\ldots,a_{k-1},(i+1)$mod $4)$ and

$((a_{0}-1)$mod $4,a_{1},\ldots,a_{k-1},(i+1)$mod $4)$ if $a_{0}$ is even.
\item[b)] $BH^{i-1}_{k}:((a_{0}+1)$mod $4,a_{1},\ldots,a_{k-1},(i-1)$mod $4)$ and

$((a_{0}-1)$mod $4,a_{1},\ldots,a_{k-1},(i-1)$mod $4)$ if $a_{0}$ is odd.
\end{enumerate}

\end{enumerate}
\end{definition}

\begin{figure}
\centering
\includegraphics[width=150mm]{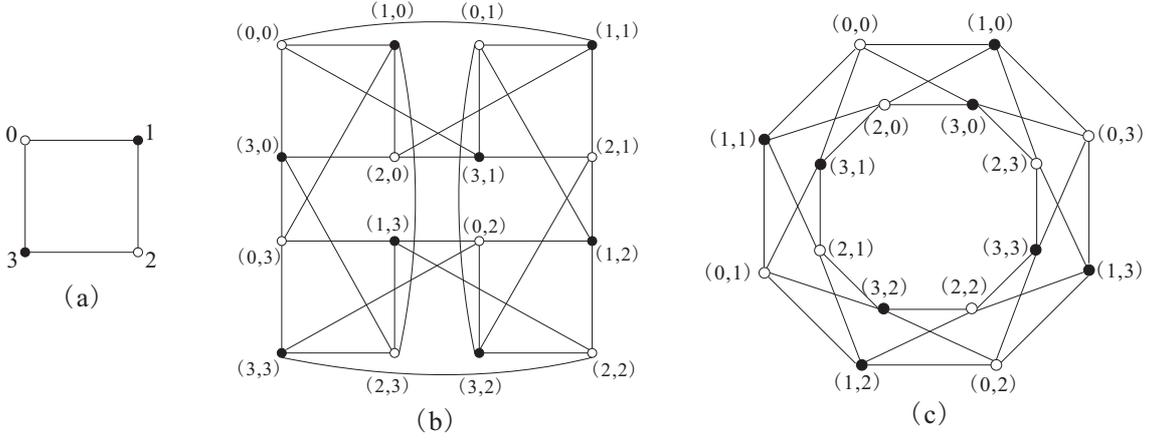}
\caption{(a) $BH_1$, (b) a layout of $BH_{2}$ and (c) ring-like layout of $BH_{2}$.} \label{g1}
\end{figure}

$BH_{1}$ is shown in Fig. \ref{g1} (a). Two distinct layouts of $BH_{2}$ are illustrated in Fig. \ref{g1} (b) and (c), respectively. Particularly, the layout of $BH_2$ in Fig. \ref{g1} (c) signifies the ring-like structure of $BH_2$.  For brevity, we shall omit ``(mod 4)'' in the rest of this paper.

Let $u$ be a neighbor of $v$ in $BH_n$. If $u$ and $v$ differ only from the inner index, then $uv$ is called a $0$-{\em dimension edge}. If $u$ and $v$ differ from $i$th outer index ($1\leq i\leq n-1$), $uv$ is called an $i$-{\em dimension edge}. It implies from Definition \ref{def1} that for each vertex $u\in V(BH_n)$, there exists two $i$-dimension neighbors, $0\leq i\leq n-1$, denoted $u^{i+}$ and $u^{i-}$, where ``+'' (resp. ``$-$'') means that the inner index of $u^{i+}$ (resp. $u^{i-}$) is that of $u$ plus one (resp. minus one). It can be deduced from Definition \ref{def2} that we can divide $BH_{n}$ into four $BH_{n-1}^{k}$s, $0\leq k \leq 3$, along dimension $n-1$. It is obvious that the edges between $BH_{n-1}^{k}$s are $(n-1)$-dimension edges. Moreover, each of $BH_{n-1}^{k}$ is isomorphic to $BH_{n-1}$. For convenience, we give some symbols as follows.

\begin{itemize}
\item $F_1$: subset of $\{\{x\}|x\in V(BH_n)\}$;
\item $F_2$: subset of $\{\{x_1,x_2\}|(x_1,x_2)\in E(BH_n)\}$;
\item $F_3$: subset of $\{\{x_1,x_2,x_3\}|(x_i,x_{i+1})\in E(BH_n)$ for each $i=1,2\}$;
\item $F_4$: subset of $\{\{x_1,x_2,x_3,x_4\}|(x_i,x_4)\in E(BH_n)$ for each $i=1,2,3\}$;
\item $Z_4$: subset of $\{\{x_1,x_2,x_3,x_4\}|(x_1,x_2),(x_2,x_3),(x_3,x_4),(x_1,x_4)\in E(BH_n)\}$.
\end{itemize}

The following basic properties of the balanced hypercube will be used in the main results of this paper.

\begin{lemma}\label{bipartite}\cite{Wu}{\bf.}
$BH_{n}$ is bipartite.
\end{lemma}

By above, vertices of odd (resp. even) inner index are colored with black (resp. white).

%\begin{lemma}\label{transitive}\cite{Wu,Zhou}{\bf.}
%$BH_{n}$ is vertex-transitive and edge-transitive.
%\end{lemma}

\begin{lemma}\label{neighbor}\cite{Wu}{\bf.}
Vertices $u=(a_{0},a_{1},\ldots,a_{n-1})$ and
$v=(a_{0}+2,a_{1},\ldots,a_{n-1})$ in $BH_{n}$ have the same neighborhood.
\end{lemma}

%
%The number of common neighbors of two vertices $u$ and $v$ is denoted by $C(u,v)$, we have the following lemma.
%
%\begin{lemma}\label{common-neighbor}\cite{Yang4}{\bf.}
%Let $u$ and $v$ be two arbitrary vertices in $BH_n$. Then $C(u,v)=0,2,$ or $2n$. Furthermore, there is exactly one vertex $w$, namely $u'$, such that $C(u,w)=2n$.
%\end{lemma}

\begin{lemma}\label{transitive}\cite{Wu,Zhou}{\bf.}
$BH_{n}$ is vertex-transitive and edge-transitive.
\end{lemma}

%\begin{lemma}\label{crossing-edges}{\rm \cite{Yang2}}{\bf.} Assume that $n\geq2$. There exist $4^{n-1}$ edges between $B^i$ and $B^{i+1}$ for each $0\leq i\leq3$.
%\end{lemma}

\begin{lemma}\label{common-neighbor}{\rm \cite{Lu2}}{\bf.}
Let $u$ and $v$ be two distinct vertices in $BH_n$. If $u$ and $v$ have a common neighbor, then $u$ and $v$ have exact two common neighbors or $2n$ common neighbors.
\end{lemma}

\begin{lemma}\label{g-connectivity-1}\cite{Yang2}{\bf.} $\kappa_0^1(BH_n)=4n-4$ for $n\geq2$.
\end{lemma}

\begin{lemma}\label{g-connectivity-2-3}\cite{Lu2}{\bf.} $\kappa_0^2(BH_n)=\kappa_0^3(BH_n)=4n-4$ for $n\geq2$.
\end{lemma}

\begin{lemma}\label{g-connectivity-4-5}\cite{D-Yang}{\bf.} $\kappa_0^4(BH_n)=\kappa_0^5(BH_n)=6n-8$ for $n\geq2$.
\end{lemma}

\section{Main results}

\subsection{$\kappa(BH_n,K_1)$, $\kappa(BH_n,K_{1,1})$, $\kappa^s(BH_n,K_1)$ and $\kappa^s(BH_n,K_{1,1})$}

It is known that $\kappa(BH_n)=2n$, so we have the following result.

\begin{theorem}\label{kappa-k1}{\bf.} $\kappa(BH_n;K_1)=2n$ and $\kappa^s(BH_n;K_1)=2n$ for $n\geq1$.
\end{theorem}

\begin{lemma}\label{kappa-k1-1-leq}{\bf.} $\kappa(BH_n;K_{1,1})\leq2n$ and $\kappa^s(BH_n;K_{1,1})\leq2n$ for $n\geq2$.
\end{lemma}
\noindent{\bf Proof.} By vertex-transitivity of $BH_n$, let $u=(0,0\cdots,0)$, $v=(1,0\cdots,0)$ and $w=(3,0\cdots,0)$. We set $F=\{v,v^{1+}\}\cup\{w,w^{1+}\}\cup\{\{u^{i+},(u^{i+})^{0+}\}|1\leq i\leq n-1\}\cup\{\{u^{i-},(u^{i-})^{0+}\}|1\leq i\leq n-1\}$. Clearly, $vv^{1+},ww^{1+}\in E(BH_n)$, and $u^{i+}(u^{i+})^{0+}, u^{i-}(u^{i-})^{0+}\in E(BH_n)$ for each $1\leq i\leq n-1$. It is obvious that $|F|=2n$. Since $N(u)\subset V(F)$, $BH_n-V(F)$ is disconnected and $u$ is one of its components. Moreover, each element in $F$ is isomorphic to $K_{1,1}$. Thus, the lemma follows.\qed

\begin{lemma}{\bf.}\label{BH-2-3-connected} If $|F_1|+|F_2|\leq 3$, then $BH_2-V(F_1\cup F_2)$ is connected.
\end{lemma}

\noindent{\bf Proof.} We may assume that $|F_1|+|F_2|=3$. Since $BH_2$ is 2-connected, $BH_2-V(F_1\cup F_2)$ is connected if $|F_2|=0$. Thus, we assume that $|F_2|\geq1$. By Lemma \ref{transitive}, we know that $BH_2$ is edge-transitive. So we assume that $u=(0,0)$, $v=(1,0)$ and $\{u,v\}\in F_2$. Let $H=BH_2-\{u,v\}$, then $H$ is 3-connected. We have the following cases.

\noindent{\bf Case 1.} $|F_2|=1$. It follows that $|F_1|=2$. Obviously, $H-V(F_1)$ is connected, which implies that $BH_2-V(F_1\cup F_2)$ is connected.

\noindent{\bf Case 2.} $|F_2|=2$. It follows that $|F_1|=1$. Pick any two adjacent vertices $x$ and $y$ in $H$, by the ring-like layout of $BH_2$, we can obtain that $H-\{x,y\}$ is 2-connected. After the deletion of any vertex in $H-\{x,y\}$, the resulting graph is connected. Thus, $BH_2-V(F_1\cup F_2)$ is connected.

\noindent{\bf Case 3.} $|F_2|=3$. We have $|F_1|=0$. By above, we know that $H-\{x,y\}$ is 2-connected. Let $x'$ and $y'$ be any two adjacent vertices in $H-\{x,y\}$. Moreover, if we delete $x'$ and $y'$ from $H-\{x,y\}$, the resulting graph is also connected. Thus, $BH_2-V(F_1\cup F_2)$ is connected.

This completes the proof. \qed

\begin{lemma}{\bf.}\label{BH-n-2n-1-connected} If $|F_1|+|F_2|\leq 2n-1$ for $n\geq2$, then $BH_n-V(F_1\cup F_2)$ is connected.
\end{lemma}

\noindent{\bf Proof.} We proceed by induction on $n$. We may assume that $|F_1|+|F_2|=2n-1$. By Lemma \ref{BH-2-3-connected}, $BH_2-V(F_1\cup F_2)$ is connected. Thus, we assume that the statement holds on $BH_i$ for each $2\leq i\leq n-1$. Next we consider $BH_n$. We set $A_k=\{x_1|\{x_1\}\in F_1\}$, $B_k^1=\{x_1|\{x_1,x_2\}\in F_2, x_1\in V(BH_{n-1}^k)$ and $x_2\not\in V(BH_{n-1}^k)\}$ and $B_k^2=\{\{x_1,x_2\}|\{x_1,x_2\}\in F_2, x_1\in V(BH_{n-1}^k)$ and $x_2\in V(BH_{n-1}^k)\}$. Clearly, $|A_k|+|B_k^1|+|B_k^2|\leq |F_1|+|F_2|=2n-1$ for each $0\leq k\leq3$. We consider the following cases.

\noindent{\bf Case 1.} $|A_k|+|B_k^1|+|B_k^2|\leq 2n-3$ for each $0\leq k\leq3$. By the induction hypothesis, each $BH_{n-1}^k-V(A_k\cup B_k^1\cup B_k^2)$ is connected. $BH_{n-1}^k$ (resp. $BH_{n-1}^{k+1}$) has $2^{2n-3}$ white (resp. black) vertices, so there are $4^{n-1}$ edges between $BH_{n-1}^{k}$ and $BH_{n-1}^{k+1}$. Since $2(2n-3)<2^{2n-3}$ whenever $n\geq3$, there exists a vertex of $BH_{n-1}^{k}$ joining to a vertex of $BH_{n-1}^{k+1}$ for each $0\leq k\leq 2$. Thus, $BH_n-V(F_1\cup F_2)$ is connected.

\noindent{\bf Case 2.} $|A_k|+|B_k^1|+|B_k^2|\geq 2n-2$ for some $0\leq k\leq3$. We may assume that $|A_0|+|B_0^1|+|B_0^2|=\max\{|A_k|+|B_k^1|+|B_k^2||0\leq k\leq3\}$, therefore, $|A_0|+|B_0^1|+|B_0^2|\geq 2n-2$. By the structure of $F_1$ and $F_2$, there may exist some $j\in\{1,3\}$ such that $|A_j|+|B_j^1|+|B_j^2|\geq2n-2$.

\noindent{\bf Case 2.1.} $|A_j|+|B_j^1|+|B_j^2|\geq2n-2$ for some $j\in\{1,3\}$. Suppose without loss of generality that $|A_1|+|B_1^1|+|B_1^2|\geq2n-2$. We claim that $|A_2|+|B_2^1|+|B_2^2|\leq1$ and $|A_3|+|B_3^1|+|B_3^2|\leq1$. Suppose not. We may assume that $|A_2|+|B_2^1|+|B_2^2|\geq2$. This implies that $|A_0|+|B_0^1|+|B_0^2|\leq2n-3$, a contradiction. Let $C$ be the subgraph induced by $\cup_{i=2}^3(V(BH_{n-1}^i)-V(F_1\cup F_2))$, then $C$ is connected. Note each black (resp. white) vertex in $BH_{n-1}^0-V(F_1\cup F_2)$ (resp. $BH_{n-1}^1-V(F_1\cup F_2)$) has a neighbor in $BH_{n-1}^3-V(F_1\cup F_2)$ (resp. $BH_{n-1}^2-V(F_1\cup F_2)$), combining the symmetry of $BH_n$, we only consider white vertices in  $BH_{n-1}^0-V(F_1\cup F_2)$. Observe that there exists a subset $F_2'\subseteq F_2$ with $|F_2'|\geq2n-3$ such that for each $\{u,v\}\in F_2'$, $u\in V(BH_{n-1}^0)$ and $v\in V(BH_{n-1}^1)$. Clearly, $u$ is a white vertex and $v$ is a black vertex. Accordingly, each white vertex of $BH_{n-1}^0-V(F_1\cup F_2)$ is adjacent to a black vertex in $BH_{n-1}^0-V(F_1\cup F_2)$, so it is connected to a vertex in $C$. Thus, $BH_n-V(F_1\cup F_2)$ is connected.

\noindent{\bf Case 2.2.} $|A_j|+|B_j^1|+|B_j^2|\leq2n-3$ for each $j\in\{1,2,3\}$. Let $C$ be the subgraph induced by $\cup_{i=1}^3(V(BH_{n-1}^i)-V(F_1\cup F_2))$, then $C$ is connected. We shall show that any vertex $u$ in $BH_{n-1}^0-V(F_1\cup F_2)$ is connected to a vertex in $C$ via a fault-free path in $BH_n$. We may assume that $u$ is a white vertex. Since $BH_n$ is triangle-free, $|N_{BH_{n-1}^0}(u)\cap V(F_1)|\leq |F_1|$ and $|N_{BH_{n-1}^0}(u)\cap V(F_2)|\leq |F_2|$. If $uu^{(n-1)+}$ or $uu^{(n-1)-}\in E(BH_n-V(F_1\cup F_2))$, we are done. Suppose not. Let $v$ be the vertex with the same neighborhood of $u$. Thus, we may assume that $\{u^{(n-1)+},(u^{(n-1)+})^{j_1+}\}\in F_2$ and $\{v,v^{(n-1)+}\}\in F_2$ for some $j_1\in\{0,\cdots,n-2\}$. Clearly, $(u^{(n-1)+})^{j_1+}\in V(BH_{n-1}^1)$. Let $D$ be the vertex set containing all $(n-1)$-dimension neighbors of  vertices in $N_{BH_{n-1}^0}(u)$. Thus, $D\subset V(BH_{n-1}^3)$ and the color of vertices in $D$ are white. Similarly, for each vertex $x\in D$, $|N_{BH_{n-1}^3}(x)\cap V(F_1)|\leq |F_1|$ and $|N_{BH_{n-1}^3}(x)\cap V(F_2)|\leq |F_2|$. Additionally, we have $|D|=|N_{BH_{n-1}^0}(u)|=2n-2>(2n-1)-2$. That is, there exists a fault-free path from $u$ to a vertex in $V(BH_{n-1}^3)-V(F_1\cup F_2)$. Thus, $BH_n-V(F_1\cup F_2)$ is connected. \qed

%\noindent{\bf Case 2.2.} $|A_0|+|B_0|=2n-1$. By the structure of $F_1$ and $F_2$, there may exist some $j\in\{1,3\}$ such that $|A_j|+|B_j^1|+|B_j^2|\geq2n-2$.
%
%\noindent{\bf Case 2.2.1.} $|A_j|+|B_j^1|+|B_j^2|\geq2n-2$ for some $j\in\{1,3\}$. Suppose without loss of generality that $|A_1|+|B_1^1|+|B_1^2|\geq2n-2$. Thus, $|A_2|+|B_2^1|+|B_2^2|<2n-3$ and $|A_3|+|B_3^1|+|B_3^2|<2n-3$. So each of $BH_{n-1}^2-V(F_1\cup F_2)$ and $BH_{n-1}^3-V(F_1\cup F_2)$ is connected, which implies that the subgraph induced by $\cup_{i=2}^3(V(BH_{n-1}^i)-V(F_1\cup F_2))$ is connected.
%
%\noindent{\bf Case 2.2.2.} $|A_j|+|B_j^1|+|B_j^2|\leq2n-3$ for each $j\in\{1,2,3\}$. Let $C$ be the subgraph induced by $\cup_{i=1}^3(V(BH_{n-1}^i)-V(F_1\cup F_2))$, then $C$ is connected. We shall show that any vertex $u$ in $BH_{n-1}^0-V(F_1\cup F_2)$ is connected to a vertex in $C$ via a fault-free path in $BH_n$.

Based on Lemmas \ref{kappa-k1-1-leq}, \ref{BH-2-3-connected} and \ref{BH-n-2n-1-connected}, we have the following theorem.

\begin{theorem}{\bf.}\label{kappa-s-k1-1} For $n\geq2$, then $\kappa^s(BH_n;K_{1,1})=2n$.
\end{theorem}

By the definitions of $\kappa(G;H)$ and $\kappa^s(G;H)$, we have $\kappa(G;H)\geq \kappa^s(G;H)$. So the following statement is straightforward.

\begin{theorem}{\bf.}\label{kappa-k1-1} For $n\geq2$, then $\kappa(BH_n;K_{1,1})=2n$.
\end{theorem}

\subsection{$\kappa(BH_n,K_{1,2})$ and $\kappa^s(BH_n,K_{1,2})$}

\begin{figure}
\centering
\includegraphics[width=72mm]{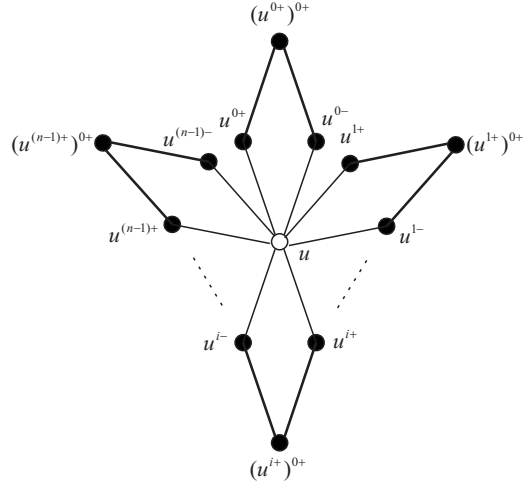}
\caption{ A $K_{1,2}$-structure-cut for $BH_n$.} \label{g2}
\end{figure}

\begin{lemma}\label{kappa-k1-2-leq}{\bf.} $\kappa(BH_n;K_{1,2})\leq n$ and $\kappa^s(BH_n;K_{1,2})\leq n$ for $n\geq2$.
\end{lemma}
\noindent{\bf Proof.} Let $u$ be an arbitrary vertex in $BH_n$. We set $F=\{\{u^{i+},(u^{i+})^{0+},u^{i-}\}|0\leq i\leq n-1\}$. We may assume that $(u^{i+})^{0+}\neq u$ if $i=0$. A $K_{1,2}$-structure-cut $F$ of $BH_n$ for $n\geq2$ is illustrated in Fig. \ref{g2}. Clearly, the subgraph induced by $u^{i+}$, $(u^{i+})^{0+}$ and $u^{i-}$ is isomorphic to $K_{1,2}$ for each $0\leq i\leq n-1$. In addition, we have $|F|=n$. Since $N(u)\subset V(F)$ and $|V(F)|=3n$, $BH_n-V(F)$ is disconnected and $u$ is one of components of $BH_n-V(F)$.X Then the lemma follows.\qed

\begin{theorem}{\bf.}\label{kappa-s-k1-2} $\kappa^s(BH_n;K_{1,2})=n$ for $n\geq2$.
\end{theorem}

\noindent{\bf Proof.} We shall show that $BH_n-V(F_1\cup F_2\cup F_3)$ is connected if $|F_1|+|F_2|+|F_3|\leq n-1$. Suppose not. Let $C$ be the smallest component of $BH_n-V(F_1\cup F_2\cup F_3)$. Note that $|V(F_1\cup F_2\cup F_2)|\leq 3n-3$, then $BH_n-V(F_1\cup F_2\cup F_3)$ is connected for $n=2$. It suffices to consider $n\geq3$. By Lemma  \ref{g-connectivity-1}, we have $4n-4>3n-3$ whenever $n\geq3$. So $|V(C)|=1$. Therefore, we assume that $x\in V(C)$. Since $BH_n$ is bipartite, $|N(x)\cap V(F_i)|\leq 2|F_i|$ for each $1\leq i\leq3$. Thus, $|N(x)\cap V(\cup _{i=1}^3 F_i)|\leq \sum_{i=1}^3|N(x)\cap F_i|\leq \sum_{i=1}^3 2|F_i|\leq 2(n-1)<2n$, which implies that there exists a neighbor of $x$ in $BH_n-V(F_1\cup F_2\cup F_3)$. So we have $|C|\geq2$, a contradiction. Thus, $BH_n-V(F_1\cup F_2\cup F_3)$ is connected. Combining $\kappa^s(BH_n;K_{1,2})\leq n$, we have $\kappa^s(BH_n;K_{1,2})=n$ for $n\geq2$. \qed

By Lemma \ref{kappa-k1-2-leq} and Theorem\ref{kappa-s-k1-2}, we have the following result.

\begin{lemma}{\bf.}\label{kappa-k1-2} $\kappa(BH_n;K_{1,2})=n$ for $n\geq2$.
\end{lemma}

\subsection{$\kappa(BH_n,K_{1,3})$ and $\kappa^s(BH_n,K_{1,3})$}

\begin{lemma}\label{kappa-k1-3-leq}{\bf.} $\kappa(BH_n;K_{1,3})\leq n$ and $\kappa^s(BH_n;K_{1,3})\leq n$ for $n\geq2$.
\end{lemma}

\begin{figure}
\centering
\includegraphics[width=80mm]{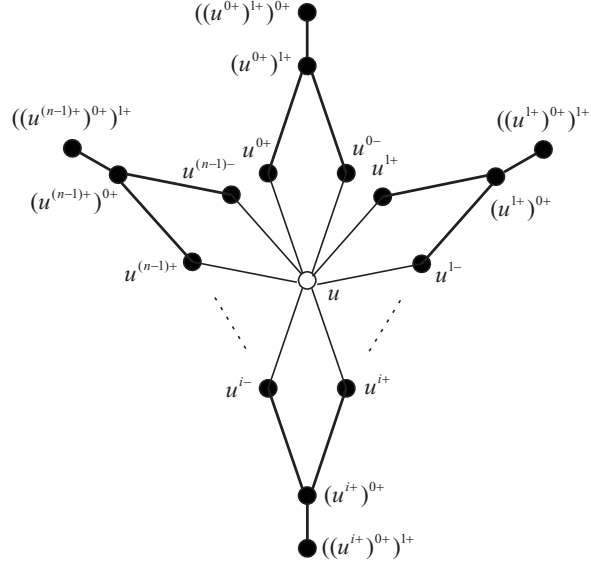}
\caption{ A $K_{1,3}$-structure-cut for $BH_n$.} \label{g3}
\end{figure}

\noindent{\bf Proof.} Let $u$ be an arbitrary vertex in $BH_n$. We set $F=\{u^{0+},(u^{0+})^{1+},((u^{0+})^{1+})^{0+}$, $u^{0-}\}\cup\{\{u^{i+},(u^{i+})^{0+},u^{i-},((u^{i+})^{0+})^{1+}\}|1\leq i\leq n-1\}$. A $K_{1,3}$-structure-cut $F$ of $BH_n$ for $n\geq2$ is shown in Fig. \ref{g3}. Clearly, the subgraph induced by $u^{0+},(u^{0+})^{1+}$, $((u^{0+})^{1+})^{0+}$ and $u^{0-}$ is isomorphic to $K_{1,3}$, and the subgraph induced by $u^{i+},(u^{i+})^{0+}$, $u^{i-}$ and $((u^{i+})^{0+})^{1+}$ is isomorphic to $K_{1,3}$ for each $1\leq i\leq n-1$. In addition, we have $|F|=n$. Since $N(u)\subset V(F)$ and $|V(F)|=4n$, $BH_n-V(F)$ is disconnected and $u$ is one of components of $BH_n-V(F)$. Then the lemma follows.\qed

\begin{theorem}{\bf.}\label{kappa-s-k1-3} $\kappa^s(BH_n;K_{1,3})=n$ for $n\geq2$.
\end{theorem}

\noindent{\bf Proof.} We shall show that $BH_n-V(\cup_{i=1}^4 F_i)$ is connected if $\sum_{i=1}^4 |F_i|\leq n-1$. Observe that $BH_2-V(\cup_{i=1}^4 F_i)$ is connected since $n-1=1$ when $n=2$. So we assume that $n\geq3$. On the contrary, suppose $BH_n-V(\cup_{i=1}^4 F_i)$ is disconnected if $\sum_{i=1}^4 |F_i|\leq n-1$. Let $C$ be the smallest component of $BH_n-V(\cup_{i=1}^4 F_i)$.

If $|F_4|\leq n-2$, then $|V(\cup_{i=1}^4 F_i)|\leq 4|F_4|+3(|F_1|+|F_2|+|F_3|)\leq 4(n-2)+3=4n-5$. By Lemma \ref{g-connectivity-1}, we have $|V(C)|=1$. Therefore, we assume that $x\in V(C)$. Since $BH_n$ is bipartite, $|N(x)\cap V(F_i)|\leq 2|F_i|$ for each $1\leq i\leq3$. We claim that there exists exact one subgraph $K_{1,3}$ of $BH_n$ such that $|N(x)\cap V(K_{1,3})|=3$. Let the center vertex of $K_{1,3}$ be $u$, and pendent vertices be $v,w$ and $y$, respectively. Accordingly, $v,w,y\in N(x)$ and $u\not\in N(x)$. Thus, $u$ and $x$ have three common neighbors, say $v,w$ and $y$. By Lemma \ref{common-neighbor}, $u$ and $x$ have $2n$ common neighbors. Therefore, $u$ and $x$ differ only the inner index. Since there exists exact one vertex $u$ such that $x$ and $u$ differ only the inner index, there exists exact one induced subgraph $K_{1,3}$ of $BH_n$ such that $|N(x)\cap V(K_{1,3})|=3$. Thus, $|N(x)\cap V(\cup _{i=1}^4 F_i)|\leq \sum_{i=1}^4|N(x)\cap V(F_i)|\leq \sum_{i=1}^3 2|F_i|+|N(x)\cap V(F_4)|\leq2+2(n-3)+3=2n-1<2n$, which implies that there exists a neighbor of $x$ in $BH_n-V(\cup_{i=1}^4 F_i)$. Thus, $BH_n-V(\cup_{i=1}^4 F_i)$ is connected.

If $|F_4|=n-1$, then $|V(F_4)|=4n-4$. By Lemmas \ref{g-connectivity-1}, \ref{g-connectivity-2-3} and \ref{g-connectivity-4-5}, $1\leq|V(C)|\leq4$. The proof of $|V(C)|=1$ is similar to that of $|F_4|\leq n-2$. Therefore, we assume that $2\leq|V(C)|\leq4$. It follows that $C$ contains at least one edge. If $2\leq|V(C)|\leq3$, combining $BH_n$ is triangle-free, it can be known that $|N(C)|>4n-4$. Note that $|V(F_4)|=4n-4$, we have a contradiction. So we assume that $|V(C)|=4$. We know that there exists at most two induced subgraphs $K_{1,3}$ of $BH_n$ such that $|N(C)\cap V(K_{1,3})|=3$ since each $K_{1,3}$ must contain a vertex in $BH_n-V(C)$ that differs only from the inner index of a vertex in $C$. We have $|N(C)\cap V(F_4)|\leq 3+3+2(n-3)<4n-4$ whenever $n\geq3$. This implies that $|V(C)|>4$, a contradiction. Thus, $BH_n-V(\cup_{i=1}^4 F_i)$ is connected.\qed

By Lemma \ref{kappa-k1-3-leq} and Theorem\ref{kappa-s-k1-3}, the following result is straightforward.

\begin{lemma}{\bf.}\label{kappa-k1-3} $\kappa(BH_n;K_{1,3})=n$ for $n\geq2$.
\end{lemma}

\subsection{$\kappa(BH_n,C_4)$ and $\kappa^s(BH_n,C_4)$}

\begin{lemma}\label{kappa-c4-leq}{\bf.} $\kappa(BH_n;C_{4})\leq n$ and $\kappa^s(BH_n;C_{4})\leq n$ for $n\geq2$.
\end{lemma}

\noindent{\bf Proof.} Let $u$ be an arbitrary vertex in $BH_n$ and let $v$ be the vertex having the same neighborhood of $u$. We set $F=\{u^{0+},(u^{0+})^{1+},u^{0-},v\}\cup\{\{u^{i+},(u^{i+})^{0+},u^{i-},(u^{i+})^{0-}\}|$ $1\leq i\leq n-1\}$. Clearly, the subgraph induced by $u^{0+},(u^{0+})^{1+},u^{0-}$ and $v$ is a 4-cycle, and the subgraph induced by $u^{i+},(u^{i+})^{0+},u^{i-}$ and $(u^{i+})^{0-}$ is also a 4-cycle for each $1\leq i\leq n-1$. In addition, we have $|F|=n$. Since $N(u)\subset V(F)$ and $|V(F)|=4n$, $BH_n-V(F)$ is disconnected and $u$ is one of components of $BH_n-V(F)$. Then the lemma follows.\qed

\begin{theorem}{\bf.}\label{kappa-s-c4} $\kappa^s(BH_n;C_{4})=n$ for $n\geq2$.
\end{theorem}

\noindent{\bf Proof.} We shall show that $BH_n-V(\cup_{i=1}^3 F_i)\cup V(Z_4)$ is connected if $\sum_{i=1}^3 |F_i|+|Z_4|\leq n-1$. Obviously, after deleting a 4-cycle or a subgraph of a 4-cycle from $BH_2$, the resulting graph is connected. So we assume that $n\geq3$. On the contrary, suppose that $BH_n-V(\cup_{i=1}^3 F_i)\cup V(Z_4)$ is disconnected if $\sum_{i=1}^3 |F_i|+|Z_4|\leq n-1$. Let $C$ be the smallest component of $BH_n-V(\cup_{i=1}^3 F_i)\cup V(Z_4)$.

If $|Z_4|\leq n-2$, then $|V(\cup_{i=1}^3 F_i)\cup V(Z_4)|\leq 4|Z_4|+3(|F_1|+|F_2|+|F_3|)\leq 4(n-2)+3=4n-5$. By Lemma \ref{g-connectivity-1}, we have $|V(C)|=1$. Thus, we assume that $x\in V(C)$. Since $BH_n$ is bipartite, $|N(x)\cap V(F_i)|\leq 2|F_i|$ for each $1\leq i\leq3$. Furthermore, $|N(x)\cap V(Z_4)|\leq 2|Z_4|$.
Therefore, $|N(x)\cap (V(\cup _{i=1}^3 F_i)\cup V(Z_4))|\leq \sum_{i=1}^3|N(x)\cap V(F_i)|+|N(x)\cap V(Z_4)|\leq 2(n-1)<2n$, which implies that there exists a neighbor of $x$ in $BH_n-V(\cup_{i=1}^3 F_i)\cup V(Z_4)$. Hence, $BH_n-V(\cup_{i=1}^3 F_i)\cup V(Z_4)$ is connected.

If $|Z_4|=n-1$, then $|V(Z_4)|=4n-4$. By Lemmas \ref{g-connectivity-1}, \ref{g-connectivity-2-3} and \ref{g-connectivity-4-5}, we have $1\leq|V(C)|\leq4$. The proof of $|V(C)|=1$ is analogous to that of $|Z_4|\leq n-2$. Therefore, we assume that $2\leq|V(C)|\leq4$. It follows that $C$ contains at least one edge, say $xy$.
If $2\leq|V(C)|\leq3$. Clearly, $|N(C)|\geq(2n-1)+(2n-2)=4n-3>|V(Z_4)|$, a contradiction. We assume that $|V(C)|=4$. If there exists one pair of vertices with the same color in $C$ differing not only the inner index, then $|N(C)|>4n-4$, a contradiction. So assume that each pair of vertices with the same color in $C$ differing only the inner index. Without loss of generality, suppose that $x,x'\in V(C)$ (resp. $y,y'\in V(C)$) and $x$ and $x'$ (resp. $y$ and $y'$) have the same neighborhood. Thus, each 4-cycle in $BH_n$ contains at most two vertices in $N(C)$, which implies that $|V(C)|>4$, a contradiction again. Thus, $BH_n-V(\cup_{i=1}^4 F_i)$ is connected.\qed

By Lemma \ref{kappa-c4-leq} and Theorem \ref{kappa-s-c4}, we have the following result.

\begin{theorem}{\bf.}\label{kappa-c4} $\kappa(BH_n;C_{4})=n$ for $n\geq2$.
\end{theorem}

\section{Conclusions}

In this paper, two novel measures of reliability and fault-tolerance, structure and substructure connectivity, are considered. For the balanced hypercube $BH_n$ ($n\geq2$), we obtain that $\kappa(BH_n;H)$ and $\kappa^s(BH_n;H)$ for $H\in\{K_1,K_{1,1},K_{1,2},K_{1,3},C_4\}$. As directions for further research, one may study $\kappa(BH_n;H)$ and $\kappa^s(BH_n;H)$ for $H\in\{P_k,C_{2k},K_{1,r}\}$ for general $k$ and $r$ with $k\geq3$ and $r\geq4$. Moreover, structure and substructure connectivity of other interconnection networks should be explored.

%\vskip 0.3 in
%
%\noindent{\bf\large Appendix}
%\vskip 0.2 in

\end{document}